\documentclass[12pt]{article} 
\usepackage{amsmath}
\usepackage{amssymb} 
\usepackage{amscd}
\usepackage{amsthm}
\usepackage[utf8]{inputenc}
\usepackage{hyperref}

\newtheorem{theorem}{Theorem}

\newtheorem{proposition}[theorem]{Proposition}
\newtheorem{lemma}[theorem]{Lemma}

\newtheorem{definition}{Definition}

\def\gr{\mathrm{gr}}

\def\dim{{\mbox{dim}}}

\def\ker{{\mbox{Ker}}}

\def\Hom {{\mbox{Hom}}}
\def\Ext {{\mbox{Ext}}}

\def\End{{\mbox{End}}}

\def\cala{{\mathcal A}}

\def\calm{{\mathcal M}}

\def\calw{{\mathcal W}}

\def\fraca{{\mathfrak A}}
 
 \def\fracg{{\mathfrak g}}

\def\bbbone{\mbox{\rm 1\hspace {-.6em} l}}

\def\Aut{{\mbox{Aut}}}

\numberwithin{equation}{section}

\begin{document}

\enlargethispage{3cm}

\thispagestyle{empty}
\begin{center}
{\bf LIE PREALGEBRAS}
\end{center}
   
\vspace{0.3cm}

\begin{center}
Michel DUBOIS-VIOLETTE
\footnote{Laboratoire de Physique Th\'eorique, UMR 8627, 
CNRS et Universit\'e Paris-Sud 11,
B\^atiment 210, F-91 405 Orsay Cedex\\
Michel.Dubois-Violette$@$u-psud.fr} and 
Giovanni LANDI \footnote{Dipartimento di Matematica e Informatica, Universit\`a di Trieste, Via A. Valerio 12/1, I-34127 Trieste, Italy and INFN, Sezione di Trieste, Trieste, Italy\\
landi$@$univ.trieste.it}\end{center}
\vspace{0,5cm}
\begin{center}
{\sl Dedicated to Henri Moscovici}
\end{center}
 \vspace{0,5cm}

%%    General info
%\subjclass[2000]{Primary: 16S37; Secondary: 16E40, 17B37}
%\keywords{Quadratic algebras. Koszul algebras. Complexes and homologies. Quantum groups. Differential calculi.}
%\date{11 June 2010; v2 27 July 2010}

%\dedicatory{Dedicated to Henri Moscovici}

\begin{abstract}
We introduce a generalization of Lie algebras within the theory of nonhomogeneous quadratic algebras and point out its relevance in the theory of quantum groups. In particular the relation between the differential calculus on quantum groups and the Koszul duality due to Positselski is made apparent.
\end{abstract}

\vfill
LPT-ORSAY 10-27
 
 \newpage
\tableofcontents

\subsection*{Acknowledgment}

GL was partially supported by the `Italian project Cofin08 - Noncommutative geometry, quantum groups and applications'.

\section{Introduction}\label{intro}

Our aim here is to describe a generalization of Lie algebras based on the theory of nonhomogeneous quadratic algebras \cite{pos:1993}, \cite{pol-pos:2005}, \cite{bra-gai:1996}, \cite{flo:2006}. 
Beside the usual Koszul duality of quadratic algebra \cite{man:1987}, \cite{man:1988}, \cite{pol-pos:2005}, a very powerful extension of it has been developed by L. Positselski \cite{pos:1993} to a duality between nonhomogeneous quadratic algebras satisfying an extension of the Poincar\'e-Birkhoff-Witt (PBW) property and curved graded differential algebras. Restricted to the quadratic-linear algebras this duality leads to graded differential algebras and this is a natural frame to understand the differential calculi on quantum groups such as the ones introduced by S.L. Woronowicz \cite{wor:1987b}, \cite{wor:1989}. In this paper we recall this theory of nonhomogeneous quadratic algebras and we specify a subclass closer to the universal enveloping algebras of Lie algebras.
This subclass of quadratic-linear algebras is characterized by a variant of the  Poincar\'e duality property for the homogeneous part refered to as the Gorenstein property \cite{art-sch:1987} and by the Koszul property.\\
In order to render the paper intelligible for readers not acquainted with the theory of nonhomogeneous quadratic algebras and to make it self-contained, an important part of this article consists of a summary of the appropriate piece of the theory of nonhomogeneous quadratic algebras.
In Section \ref{Lie}, we recall some basic facts on Lie algebras, associated complexes and the PBW property. Section \ref{aalg} provides a similar discussion for associative algebras instead of Lie algebras. Quadratic algebras, Koszul complexes, etc. are reviewed in Section \ref{quad}. Section \ref{ExI} introduces some representative examples of quadratic algebras which will be used later as the homogeneous part of examples of nonhomogeneous quadratic algebras. In Section \ref{Nonh} the nonhomogeneous quadratic algebras are introduced and the PBW-property is defined for these algebras. Section \ref{dual} is devoted to the Koszul duality of Positselski. In Section \ref{ExII} nonhomogeneous versions of the examples of Section \ref{ExI} are discussed. Lie prealgebras are defined in Section \ref{pseudoL} and are shown to be in duality with the differential quadratic Koszul Frobenius algebras. In Section \ref{coh} the corresponding generalization of Chevalley-Eilenberg complexes is investigated.\\
Throughout the paper, $\mathbb K$ is a (commutative) field and all vector spaces and algebras are over $\mathbb K$. By an algebra without other mention, we always mean an associative algebra and by a graded algebra we mean a $\mathbb N$-graded algebra.

\section{Lie algebras, PBW property, etc.}\label{Lie}

Let $E$ be a finite-dimensional vector space equipped with an antisymmetric bilinear product
\begin{equation}
(x,y)\mapsto [x,y]
\label{brac}
\end{equation}
that is one has a linear mapping 
\[
\psi:\wedge^2 E\rightarrow E
\]
and the product is given by
\[
[x,y]+\psi(x\wedge y)=0
\]
for any $x,y\in E$.\\

Let the tensor algebra $T(E)$ over $E$ be equipped with the filtration
\[
F^n(T(E))=\otimes_{m\leq n} E^{\otimes^m}
\]
associated to its graduation and let $\fraca$ be the algebra defined by 
\[
\fraca=T(E)/(\{x\otimes y-y\otimes x-[x,y]\vert x,y\in E\})
\]
where $(P)$ denotes for any $P\subset T(E)$ the two-sided ideal of $T(E)$ generated by $P$. 
Notice that one has
\[
(\{x\otimes y-y\otimes x -[x,y]\vert x,y\in E\})=(\{r+\psi(r)\vert r\in \wedge^2 E\})
\]
where in the right-hand side of this equality $\wedge^2 E$ is considered as a subset of $E\otimes E$.\\

The filtration of the tensor algebra induces a filtration $F^n(\fraca)$ of $\fraca$ and one defines the associated graded algebra
\[
\gr (\fraca)=\oplus_n F^n(\fraca)/F^{n-1}(\fraca)
\]
(with the convention $F^n(\fraca)=0$ whenever $n<0$). One has a canonical surjective homomorphism
\begin{equation}
\text{can}:S(E)\rightarrow \gr (\fraca)
\label{canS}
\end{equation}
of graded algebras of the symmetric algebra $S(E)$ onto $\gr(\fraca)$.\\

Consider the transpose $\psi^t$ of $\psi$, that is the linear mapping
\[
\psi^t:E^\ast\rightarrow \wedge^2E^\ast
\]
where $E^\ast$ denotes the dual vector space of $E$. This transpose $\psi^t$ has a unique extention as an antiderivation
\begin{equation}
d:\wedge E^\ast\rightarrow \wedge E^\ast
\label{antiwedge}
\end{equation}
of the exterior algebra $\wedge E^\ast$ over $E^\ast$. One has the following theorem
\begin{theorem}\label{LIE}
The following conditions $\mathrm{(a), (b)}$ and $\mathrm{(c)}$ are equivalent :\\

\noindent $\mathrm{(a)}$ the bracket $(\ref{brac})$ satisfies the Jacobi identity, i.e. one has
\[
[[x,y],z]+[[y,z],x]+[[z,x],y]=0
\]
for any $x,y,z\in E$,\\

\noindent $\mathrm{(b)}$ the antiderivation $(\ref{antiwedge})$ is a differential of $\wedge E^\ast$ i.e. $d^2=0$,\\

\noindent $\mathrm{(c)}$ the canonical homomorphism $(\ref{canS})$ is an isomorphism of $S(E)$ onto $\gr(\fraca)$.
\end{theorem}
This theorem is more or less classical and easy to prove (see for instance \cite{mdv:2001}).\\

Condition (a) means that $(E,[\bullet,\bullet])$ is a Lie algebra, Condition (b) means that $(\wedge E^\ast, d)$ is a graded differential algebra while Condition (c) means that one has the  Poincar\'e-Birkhoff-Witt (PBW) property.\\

In fact $\fraca$ is then the universal enveloping algebra $U(E)$ of the Lie algebra $E$ while the graded differential algebra $\wedge E^\ast$ is the basic complex to construct the Chevalley-Eilenberg complexes ($d$ is the Koszul differential).

\section{Associative algebras}\label{aalg}

Assume now that the finite-dimensional vector space $E$ is equipped with an arbitrary bilinear product
\begin{equation}
(x,y)\mapsto xy
\label{prod}
\end{equation}
so one has a linear mapping
\begin{equation}
\psi:E\otimes E\rightarrow E
\label{psiA}
\end{equation}
and the product is given by
\[
xy+\psi(x\otimes y)=0
\]
for any $x,y\in E$. Thus $E$ is an algebra which may be nonassociative.\\

Let then $\fraca$ be the (associative) algebra
\[
\fraca=T(E)/(\{x\otimes y-xy\vert x,y\in E\})
\]
with the same conventions as before. One has
\[
(\{x\otimes y-xy\vert x,y\in E\})=(\{r+\psi(r)\vert r\in E\otimes E\})
\]
for the defining ideal of $T(E)$.\\

The algebra $\fraca$ is filtered by the filtration induced from the one of the tensor algebra $T(E)$ and one has the associated graded algebra.
\[
\gr(\fraca)=\otimes_n F^n(\fraca)/F^{n-1}(\fraca)
\]
as in last section. Let us investigate the structures of $\fraca$ and $\gr(\fraca)$. The projection of $T(E)$ onto the degree 0 induces homomorphisms 
\[
\varepsilon : T(E)\rightarrow \mathbb K
\]
and 
\[
\varepsilon : \fraca\rightarrow\mathbb K
\]
of unital (associative) algebras. The kernels $T_+(E)$ and $\fraca_+$ of these homomorphisms are two-sided ideals and one has the decompositions
\[
T(E)=T_+(E)\oplus \mathbb K\bbbone
\]
and
\[
\fraca=\fraca_+\oplus \mathbb K\bbbone
\]
of these algebras. By construction and by the very universal property of $T_+(E)$, one has a surjective homomorphism of algebra
\[
\alpha:E \rightarrow \fraca_+
\]
which is such that any algebra-homomorphism of $E$ into an associative algebra $A$ factorizes through a unique homomorphism of $\fraca_+$ into $A$ and $\alpha$. It follows that $\fraca_+$ is the quotient of $E$ by the two-sided ideal $I$ generated by the associators
\[
[x,y,z]=(xy)z-x(yz)
\]
for $x,y,z\in E$. One has
\[
F^n(\fraca)=\fraca\ \text{for}\ n\geq 1\ \text{and}\ F^0(\fraca)=\mathbb K\bbbone
\]
so the graded algebra $\gr(\fraca)$ is given by
\[
\gr(\fraca)=\gr_1(\fraca)\oplus \mathbb K \bbbone
\]
where $\gr_1(\fraca)=\fraca_+=E/I$ as vector space but where the product of 2 elements of $\gr_1(\fraca)$ vanishes. Let $\cala$ be the graded algebra
\[
\cala=T(E)/(E\otimes E)=E\oplus \mathbb K\bbbone
\]
with obvious notations. It is clear that one has a canonical surjective homomorphism
\begin{equation}
\text{can}: \cala\rightarrow \gr(\fraca)
\label{canA}
\end{equation}
of graded algebra which is an isomorphism if and only if $E$ is an associative algebra.\\

Consider the transpose linear mapping
\[
\psi^t:E^\ast\rightarrow E^\ast \otimes E^\ast
\]
of $\psi$. This transpose $\psi^t$ has a unique extension as an antiderivation
\begin{equation}
d:T(E^\ast)\rightarrow T(E^\ast)
\label{antiT}
\end{equation}
of the tensor algebra over $E^\ast$. One has the following analog of Theorem \ref{LIE}.

\begin{theorem}\label{ALG}
The following conditions $\mathrm{(a)}$, $\mathrm{(b)}$ and $\mathrm{(c)}$ are equivalent: \\

\noindent $\mathrm{(a)}$ the product $\mathrm{(\ref{prod})}$ is associative,\\

\noindent $\mathrm{(b)}$ the antiderivation $\mathrm{(\ref{antiT})}$ is a differential i.e. $d^2=0$,\\

\noindent $\mathrm{(c)}$ the canonical homomorphism $\mathrm{(\ref{canA})}$ is an isomorphism.
\end{theorem}

We have explained above the equivalence $(a)\Leftrightarrow (c)$. For the equivalence $(a)\Leftrightarrow (b)$ as in Theorem \ref{LIE} see for instance in 
\cite{mdv:2001}.\\

In fact $\fraca$ is then the algebra
\[
\fraca=E\oplus \mathbb K\bbbone=\tilde E
\]
that is the algebra obtained by adjoining a unit element to $E$; the algebra $E$ itself identifies with the two-sided ideal $\fraca_+$ of $\fraca$, i.e. one has $E=\fraca_+$.
Notice that then the graded differential algebra $(T(E^\ast), d)$, or more precisely $(T_+(E^\ast),d)$, is the basic complex to construct the Hochschild complexes.\\

We shall see that the algebras $\fraca$ of Section \ref{Lie} satisfying the conditions of Theorem \ref{LIE} and  of this section satisfying the conditions of Theorem \ref{ALG} are both examples of Koszul quadratic-linear algebras but that the latter one does not fall into the class of Lie prealgebras.

\section{Quadratic algebras}\label{quad}

A (homogeneous) {\sl quadratic algebra} \cite{man:1987}, \cite{man:1988}, \cite{pol-pos:2005} is an associative algebra $\cala$ of the form
\[
\cala=A(E,R)=T(E)/(R)
\]
where $E$ is a finite-dimensional vector space, $R$ is a subspace of $E\otimes E$ and where $(R)$ denotes the two-sided ideal of the tensor algebra $T(E)$ over $E$ generated by $R(\subset E\otimes E)$. Space $E$ is the space of generators of $\cala$ and subspace $R$ of $E\otimes E$ is the space of relations of $\cala$. The algebra $\cala=A(E,R)$ is naturally a graded algebra $\cala=\oplus_{n\in \mathbb N}\cala_n$ which is connected, i.e. such that $\cala_0=\mathbb K\bbbone$ and generated in degree 1, $\cala_1=E$.\\
To a quadratic algebra $\cala=A(E,R)$ as above one associates another quadratic algebra, {\sl its Koszul dual} $\cala^!$, defined by
\[
\cala^!=A(E^\ast,R^\perp)
\]
where $E^\ast$ denotes the dual vector space of $E$ and $R^\perp\subset E^\ast \otimes E^\ast$ is the orthogonal of $R\subset E\otimes E$ defined by 
\[
R^\perp = \{\omega \in E^\ast\otimes E^\ast\vert \langle\omega,r\rangle=0,\forall r\in R\}
\]
where, by using the finite-dimensionality of $E$, we have identified $E^\ast \otimes E^\ast$ with the dual vector space $(E\otimes E)^\ast$ of $E\otimes E$. One has of course
\[
(\cala^!)^!=\cala
\]
and the dual vector spaces $\cala^{!\ast}_n$ of the homogeneous components $\cala^!_n$ of $\cala^!$ are given by
\[
\cala^{!\ast}_1=E
\]
and
\begin{equation}
\cala^{!\ast}_n=\cap_{r+s+2=n} E^{\otimes^r}\otimes R \otimes E^{\otimes^s}
\label{En}
\end{equation}
for $n\geq 2$ as easily verified. In particular one has $\cala^{!\ast}_2=R$ and $\cala^{!\ast}_n\subset E^{\otimes^n}$ for any $n\in \mathbb N$.\\
Consider the sequence of free left $\cala$-modules
\begin{equation}
\stackrel{b}{\rightarrow}\cala\otimes \cala^{!\ast}_{n+1}\stackrel{b}{\rightarrow}\cala\otimes \cala^{!\ast}_n\rightarrow \dots \rightarrow \cala\otimes \cala^{!\ast}_2\stackrel{b}{\rightarrow}\cala\otimes E\stackrel{b}{\rightarrow}\cala\rightarrow 0
\label{K}
\end{equation}
where $b:\cala\otimes \cala^{!\ast}_{n+1}\rightarrow \cala\otimes \cala^{!\ast}_n$ is induced by the left $\cala$-module-homomorphism of $\cala\otimes E^{\otimes^{n+1}}$ into $\cala\otimes E^{\otimes^n}$ defined by
\[
b(a\otimes (x_0\otimes x_1\otimes \dots \otimes x_n))=ax_0 \otimes (x_1\otimes \dots \otimes x_n)
\]
for $a\in \cala$, $x_i\in E$. It follows from (\ref{En}) that $\cala^{!\ast}_n\subset R\otimes E^{\otimes^{n-2}}$ for $n\geq 2$ which implies that $b^2=0$. As a consequence the sequence (\ref{K}) is a chain complex of free left $\cala$-modules called the {\sl Koszul complex} of the quadratic algebra $\cala$ and denoted by $K(\cala)$. The quadratic algebra $\cala$ is said to be a {\sl Koszul algebra} whenever its Koszul complex is acyclic in positive degrees, i.e. whenever $H_n(K(\cala))=0$ for $n\geq 1$. One shows easily that $\cala$ is a Koszul algebra if and only if its Koszul dual $\cala^!$ is a Koszul algebra.\\
It is important to realize that the presentation of $\cala$ by generators and relations is equivalent to the exactness of the sequence
\begin{equation}
\cala\otimes R \stackrel{b}{\rightarrow} \cala\otimes E \stackrel{b}{\rightarrow}\cala \stackrel{\varepsilon}{\rightarrow} \mathbb K\rightarrow 0
\label{Pres}
\end{equation}
so one always has
\[
H_1(K(\cala))=0\ \text{and}\ H_0(K(\cala))=\mathbb K
\]
and 
\begin{equation}
K(\cala)\stackrel{\varepsilon}{\rightarrow} \mathbb K \rightarrow 0
\label{KRes}
\end{equation}
is a free resolution of the trivial module $\mathbb K$ whenever $\cala$ is Koszul.This resolution is then a minimal projective resolution of $\mathbb K$ in the graded category (i.e. the category of graded modules) \cite{car:1958}. In the above sequences, $\varepsilon$ is induced by the projection onto degree 0.\\

Let $\cala=A(E,R)$ be a quadratic Koszul algebra such that one has $\cala^!_D\not=0$ and $\cala^!_n=0$ for $n>D$. It follows that the trivial (left) module $\mathbb K$ has projective dimension $D$ which implies that $\cala$ has global dimension $D$ (see \cite{car:1958}). It is worth noticing that this also implies that the Hochschild dimension of $\cala$ is $D$ (see \cite{ber:2005}). By applying the functor $\Hom_\cala(\bullet, \cala)$ to the Koszul chain complex $K(\cala)$ of left $\cala$-modules one obtains the cochain complex $L(\cala)$ of right $\cala$-modules 
\begin{equation}
0\rightarrow \cala\stackrel{b'}{\rightarrow} \dots\stackrel{b'}{\rightarrow}\cala^!_n \otimes \cala \stackrel{b'}{\rightarrow} \cala^!_{n+1} \otimes \cala \stackrel{b'}{\rightarrow} 
\label{L}
\end{equation}
$b'$ being the left multiplication by $\sum_k\theta^k\otimes e_k$ in $\cala^!\otimes \cala$ where $(e_k)$ is a basis of $E$ with dual basis $(\theta^k)$. The algebra $\cala$ is said to be {\sl Koszul-Gorenstein} if it is Koszul of finite global dimension $D$ as above and if $H^n(L(\cala))=\mathbb K\delta^n_D$. Notice that this implies that $\cala^!_n\simeq \cala^{!\ast}_{D-n}$ as vector spaces (a version of the  Poincar\'e duality).

\section{Examples I}\label{ExI}

In this section we analyse some examples of quadratic algebras $\cala=A(E,R)$ which will appear as the ``homogeneous parts" of representative examples of the concepts described in the sequel. 
\\

\subsection*{ 1. Case $R=\wedge^2E\subset E\otimes E$ } ~\\
In the case where $R=\wedge^2E$ considered as a subspace of $E\otimes E$, the quadratic algebra $\cala=A(E,R)$ is the symmetric algebra $S(E)$ over $E$
\[
\cala=A(E,\wedge^2 E)=S(E)
\]
and one has
\[
\cala^!=\wedge E^\ast
\]
for its Koszul dual.\\

In this case the Koszul complex is $(S(E)\otimes \wedge^\bullet E,b)$ and the exterior differential $d$ is such that
\[
db + bd =\ \text{degree}\ (\text{in}\ S^\bullet\ \text{and}\ \wedge^\bullet)
\]
so $d$ gives an homotopy for $b$ in positive degrees. Therefore $\cala=S(E)$ is Koszul and, since $\wedge^DE\not=0$ and $\wedge^nE=0$ for $n>D=\dim (E)$, $\cala=S(E)$ has global dimension $D=\dim(E)$. It is not hard to show that $\cala=S(E)$ is in fact Koszul-Gorenstein.\\

The Koszul dual $\cala^!=\wedge E^\ast$ is therefore also Koszul but is not of finite global dimension since its Koszul complex $(\wedge E^\ast \otimes S^\bullet(E), b)$ has components of degree $n$
\[
\wedge E^\ast \otimes S^n(E)\not=0
\]
for any $n\in \mathbb N$.
\\

\subsection*{ 2. Case $R=E\otimes E$} ~\\
In the case where $R=E\otimes E$, $\cala=A(E,R)$ is the trivial quadratic algebra
\[
\cala=A(E,E\otimes E)=E\oplus \mathbb K\bbbone
\]
with $\cala_1=E$ and $\cala_n=0$ for $n\geq 2$, and one has
\[
\cala^!=A(E^\ast,\{0\})=T(E^\ast)
\]
for its Koszul dual, i.e. it is the tensor algebra over $E^\ast$.\\

The Koszul complex for $\cala^!=T(E^\ast)$ reduces then to
\[
0\rightarrow T(E^\ast)\otimes E^\ast \stackrel{b}{\rightarrow} T(E^\ast)\rightarrow 0
\]
where $b$ is the ``product" so $\cala^!=T(E^\ast)$ is Koszul of global dimension 1 but is obviously not Koszul-Gorenstein whenever $\dim(E)\geq 2$.\\

The quadratic algebra $\cala=E\oplus \mathbb K\bbbone$ is thus also Koszul but, since its Koszul complex $((E\oplus\mathbb K\bbbone)\otimes T^\bullet(E),b)$ has components of degree $n$

\[
(E\oplus \mathbb K \bbbone) \otimes E^{\otimes^n} \not= 0
\]
for any $n\in \mathbb N$, it is not of finite global dimension.
\\

\subsection*{ 3. Example connected with the calculus  \cite{wor:1987b} on the twisted $SU(2)$ group} ~\\
Let $\cala$ be the quadratic algebra generated by 3 elements $\nabla_0,\nabla_1,\nabla_2$ of degree 1 with relations

\begin{equation}
\left\{
\begin{array}{l}
\nu\nabla_2\nabla_0-\displaystyle{\frac{1}{\nu}}\nabla_0\nabla_2 = 0\\
\\
\nu^2\nabla_1\nabla_0-\displaystyle{\frac{1}{\nu^2}}\nabla_0\nabla_1 = 0\\
\\
\nu^2\nabla_2\nabla_1-\displaystyle{\frac{1}{\nu^2}}\nabla_1\nabla_2 = 0
\end{array}
\right.
\label{rWh}
\end{equation}
where $\nu\in \mathbb K\backslash \{0\}$. By setting 
\begin{equation}
E=\mathbb K\nabla_0 \oplus \mathbb K \nabla_1\oplus \mathbb K \nabla_2
\label{EWh}
\end{equation}
and

\begin{eqnarray}
R & = & \mathbb K(\nu \nabla_2\otimes \nabla_0-\frac{1}{\nu}\nabla_0 \otimes \nabla_2) \oplus  \mathbb K(\nu^2\nabla_1\otimes \nabla_0-\frac{1}{\nu^2}\nabla_0\otimes \nabla_1)\nonumber\\
& \oplus & \mathbb K(\nu^2\nabla_2 \otimes \nabla_1	-\frac{1}{\nu^2}\nabla_1\otimes \nabla_2)
\label{RWh}
\end{eqnarray}

one has $\cala=A(E,R)$.\\

The Koszul dual of $\cala$ is the algebra $\cala^!$ generated by 3 elements $\omega_0,\omega_1,\omega_2$ of degree 1 with relations
\begin{equation}
\left\{
\begin{array}{l}
(\omega_\alpha)^2=0,\>\>  \forall\alpha\in \{0,1,2\}\\
\omega_2\omega_0 + \nu^2 \omega_0\omega_2 =0\\
\omega_1\omega_0+\nu^4\omega_0\omega_1=0\\
\omega_2\omega_1+\nu^4\omega_1\omega_2=0
\end{array}
\right.
\label{KdrWR}
\end{equation}
The algebra $\cala$ is a deformation of $S(E)$ with $\dim(E)=3$, it is Koszul 
\cite{gur:1990} of global dimension 3 and is Koszul-Gorenstein (see however  Remark c) below). The argument for the Koszul property is similar to the one used in the example 1 above (see also \cite{wam:1993}).\\
The Koszul algebra $\cala^!$ has infinite global dimension.
\\

\subsection*{ 4. Generalization : $q$-deformed polynomial algebras}~\\
Let $\cala$ be the algebra generated by $D$ elements $X^\lambda$ ($\lambda\in \{1,\dots,D\}$) with relations
\begin{equation}
X^\mu X^\nu-q^{\mu\nu} X^\nu X^\mu=0
\label{qSD}
\end{equation}
for $\mu,\nu\in \{1,\dots,D\}$ where the $q^{\mu\nu}\in \mathbb K$ are such that 
\begin{equation}
q^{\mu\nu} q^{\nu\mu}=1,\ \ q^{\lambda\lambda}=1
\label{qrs}
\end{equation}
for $\lambda, \mu,\nu\in \{1,\dots, D\}$. Again this quadratic algebra is Koszul 
\cite{gur:1990} of global dimension $D$ and is in fact Koszul-Gorenstein. Its Koszul dual $\cala^!$ is generated by the $D$ elements $\omega_\alpha$ ($\alpha\in \{1,\dots, D\}$) with relations
\begin{equation}
\omega_\mu\omega_\nu + q^{\nu\mu} \omega_\nu \omega_\mu=0
\label{qwedgeD}
\end{equation}
for $\mu,\nu\in \{1,\dots, D\}$. It is again a Koszul algebra of infinite global dimension.\\

\noindent\underbar{Remarks}\\

\noindent a) Let $\cala$ be a quadratic Koszul algebra with Koszul dual $\cala^!$. 
Then their  Poincar\'e series satisfy 
\[
P_\cala(t)P_{\cala^!}(-t)=1
\]
where the {\sl  Poincar\'e series} $P_\cala(t)$ is defined by 
\[
P_\cala(t)=\sum_n \dim(\cala_n)t^n
\]
for any graded algebra $\cala=\oplus_n \cala_n$.\\
It follows that if $\cala$ has a finite global dimension $D\geq 1$ then $\cala^!$ is of infinite global dimension.\\

\noindent b) Recall that a graded algebra $\cala=\oplus_n \cala_n$ is said to have {\sl polynomial growth} whenever there are positive $K$ and $N\in \mathbb N$ such that
\[
\dim (\cala_n)\leq K n^{N-1}
\]
for any $n\in \mathbb N$. The algebras $\cala$ of Examples 1, 2, 3 and 4 above have polynomial growth but the tensor algebra $T(E)$ has exponential growth whenever $\dim(E)\geq 2$.\\

\noindent c) To complete the proof of the Koszul and the Koszul-Gorenstein properties for the algebras of Examples 3 and 4 above, one must use in addition to the results of \cite{gur:1990} the fact that the Koszul and the Gorenstein properties are stable by the twists \cite{art-tat-vdb:1991}. This was proved in \cite{pot:2006} in the more general context of homogeneous algebras.

\section{Nonhomogeneous quadratic algebras }\label{Nonh}

In this subsection we let $E$ be as before a finite-dimensional vector space and we endow the tensor algebra $T(E)$ with its natural filtration $F^n(T(E))=\oplus_{m\leq n} E^{\otimes^m}$, $(n\in \mathbb N)$.\\
A {\sl nonhomogeneous quadratic algebra} \cite{pos:1993}, \cite{bra-gai:1996}, \cite{flo:2006} is an algebra $\fraca$ of the form
\[
\fraca=A(E,P)=T(E)/(P)
\]
where $P$ is a subspace of $F^2(T(E))$ and where $(P)$ denotes as above the two-sided ideal of $T(E)$
generated by $P$. The filtration of the tensor algebra $T(E)$ induces a filtration $F^n(\fraca)$ of $\fraca$ and one associates as before to $\fraca$ the graded algebra
\[
\gr (\fraca)=\oplus_n F^n(\fraca)/F^{n-1}(\fraca)
\]
which is refered to as the {\sl associated graded algebra} to the filtered algebra $\fraca$. Let $R$ be the image of $P$ under the canonical projection of $F^2(T(E))$ onto $E\otimes E$ and let $\cala=A(E,R)$ be the homogeneous quadratic algebra $T(E)/(R)$; $\cala$ will be refered to as the {\sl quadratic part} of $\fraca$. There is again the canonical  surjective homomorphism 
\[
\text{can} :\cala\rightarrow \gr(\fraca)
\]
 of graded algebras and $\fraca$ is said to have the {\sl  Poincar\'e-Birkhoff-Witt (PBW) property} whenever this homomorphism is an isomorphism. This terminology  comes from the example where $\fraca$ is the universal enveloping algebra $U(\fracg)$ of a Lie algebra $\fracg$ (see in  Section \ref{Lie}).\\

One has the following theorem  \cite{bra-gai:1996}.

\begin{theorem} \label{BG}
Let $\fraca$ and $\cala$ be as above. If $\fraca$ has the PBW property then the following conditions $\mathrm{(i)}$ and $\mathrm{(ii)}$ are satisfied :\\
\\ $\mathrm{(i)}$ $P\cap F^1 (T(E))=0$,\\

\noindent $\mathrm{(ii)}$ $(P.E+E.P)\cap F^2 (T(E)) \subset P$.\\

\noindent Assume that $\cala$ is a Koszul algebra, then conversely if conditions $\mathrm{(i)}$ and $\mathrm{(ii)}$ are satisfied $\fraca$ has the PBW property.
\end{theorem}
Condition (i) means that $P$ is obtained from $R$ by adding to each non zero element of $R$ terms of degrees 1 and 0. In other words, Condition (i) means that there are linear mappings $\psi_1:R\rightarrow E$ and $\psi_0:R\rightarrow \mathbb K$ such that one has
\begin{equation}
P=\{x+\psi_1(x)+\psi_0(x)\bbbone \vert x\in R\}
\label{Ci}
\end{equation}
which gives $P$ in terms of $R$. Condition (ii) which is a generalization of the Jacobi identity is then given more explicitely by the following proposition (see e.g. in \cite{pol-pos:2005}).\\

\noindent\underbar{Remark}. In the second part of Theorem \ref{BG} one can replace the condition of Koszulity of $\cala$ by a weaker condition, namely the condition of acyclicity in degree 3 and 2 of the Koszul complex of $\cala$ to conclude that $\fraca$ has the PBW property whenever conditions (i) and (ii) are satisfied.\\

In spite of the above remark, it turns out that the Koszul property for $\cala$ is a very natural and helpful property. A nonhomogeneous quadratic algebra $\fraca$ is said to be {\sl Koszul} whenever it has the PBW property and its (homogeneous) quadratic part $\cala$ is Koszul.

\begin{proposition} \label{EPBW}

Assume that Condition $\mathrm{(i)}$ of Theorem \ref{BG} is satisfied, let $\psi_1$ and $\psi_0$ be as above and set $\calw_3=(R\otimes E)\cap (E\otimes R)$. Then Condition $\mathrm{(ii)}$ of Theorem \ref{BG} is
equivalent to the following conditions $\mathrm{(a)}$, $\mathrm{(b)}$ and $\mathrm{(c)}$\\

\noindent $\mathrm{(a)}$ $(\psi_1\otimes I - I\otimes \psi_1)(\calw_3)\subset R$\\

\noindent $\mathrm{(b)}$ $\big(\psi_1\circ (\psi_1\otimes I-I\otimes \psi_1)- (\psi_0 \otimes I-I\otimes \psi_0)\big)(\calw_3)=0$\\

\noindent $\mathrm{(c)}$ $\psi_0\circ (\psi_1 \otimes I-I\otimes \psi_1)(\calw_3)=0$\\

\noindent  where $I$ is the identity mapping of $E$ onto itself.
\end{proposition}

\section{Duality}\label{dual}

In this section we recall the Koszul duality of \cite{pos:1993} for the nonhomogeneous quadratic Koszul algebras.\\

Let $\fraca=T(E)/(P)$ be a nonhomogeneous quadratic algebra and let its quadratic part  be 
$\cala=T(E)/(R)$. We assume that the condition (i) of Theorem \ref{BG} is satisfied and we let $\psi_1:R\rightarrow E$ and $\psi_0:R\rightarrow \mathbb K$ be the corresponding linear mappings as in (\ref{Ci}).\\

Let $\psi^t_1:E^\ast\rightarrow R^\ast$ and $\psi^t_0:\mathbb K^\ast \rightarrow R^\ast$ be the transpose linear mappings of $\psi_1$ and $\psi_2$. The Koszul dual $\cala^!=\oplus_n\cala^!_n$ of $\cala$ is such that 
\[
\cala^!_n=\cap_{r+s+2=n} (E^{\otimes^r}\otimes R\otimes E^{\otimes^s})^\ast
\]
for $n\geq 2$ and $\cala^!_1=E^\ast$ so in particular we can write 
\begin{equation}
\psi^t_1:\cala^!_1\rightarrow \cala^!_2
\label{tpsi1}
\end{equation}
and 
\begin{equation}
\psi^t_0(1)=F
\label{curvat}
\end{equation}
is an element of  $\cala^!_2$. One has $\cala^!_3=\calw^\ast_3$ and Condition (a) of Proposition \ref{EPBW} means that $\psi^t_1$ extends as an antiderivation 
\begin{equation}
d:\cala^!\rightarrow \cala^!
\label{antiK}
\end{equation}
 of degree 1 of $\cala^!$, Condition (b) reads then
\begin{equation}
d^2 \alpha = [F,\alpha]
\label{Curvd}
\end{equation}
for any $\alpha\in \cala^!$ while Condition (c) reads
\begin{equation}
dF=0
\label{Bianchi}
\end{equation}
for the element $F=\psi_0^t (1)$ of $\cala^!_2$. \\

A graded algebra equipped with an antiderivation $d$ of degree 1 and an element $F$ of degree 2 satisfying the conditions (\ref{Curvd}) and (\ref{Bianchi}) above is called a {\sl curved graded differential algebra} \cite{pos:1993}. The above correspondence between nonhomogeneous quadratic algebras satisfying conditions (i) and (ii) of Theorem \ref{BG} and curved graded differential algebras is an anti-equivalence of categories between the category of nonhomogeneous quadratic algebras satisfying conditions (i) and (ii) and a full subcategory of the category of curved graded differential algebras \cite{pol-pos:2005}.\\

\noindent Let us introduce a more specific setting.\\
A (curved) graded differential algebra which is as graded algebra a quadratic algebra endowed with its natural graduation will be refered to as a {\sl (curved) differential quadratic algebra}. If furthermore the underlying quadratic algebra is Koszul then it will be refered to as a {\sl curved differential quadratic Koszul algebra}.\\
Given two curved differential quadratic algebras, $\cala=(A(E,R),d,F)$ and $\cala'=(A(E',R'),d',F')$ a {\sl morphism of $\cala$ to $\cala'$} is a linear mapping 
\[
\alpha:E\rightarrow E'
\]
such that
\[
(\alpha \otimes\alpha) (R)\subset R'
\]
and such that the corresponding unital algebra homomorphism $\tilde\alpha:\cala\rightarrow \cala'$ satisfies
\[
\tilde\alpha(F)=F'
\]
and
\[
\tilde\alpha(dx)=d'\tilde\alpha(x)
\]
for any $x\in \cala$. One has then the corresponding category of curved differential quadratic algebras and the full subcategory of curved differential quadratic Koszul algebras.\\
Given two nonhomogeneous quadratic algebras $\fraca=A(E,P)=T(E)/(P)$ and $\fraca'=A(E',P')=T(E')/(P')$ a {\sl morphism from $\fraca$ to $\fraca'$}
\[
\alpha:\fraca\rightarrow \fraca'
\]
is a linear mapping 
\[
\alpha:E\rightarrow E'
\]
such that
\[
T(\alpha)(P)\subset P'
\]
where $T(\alpha):T(E)\rightarrow T(E')$ is the corresponding unital algebra homomorphism. One has the corresponding category of nonhomogeneous quadratic algebras and the full subcategory of nonhomogeneous quadratic Koszul algebras.\\

It is clear that the correspondence $\fraca\mapsto (\cala^!,d,F)$ defines a contravariant functor from the category of nonhomogeneous quadratic algebras satisfying the condition (i) and (ii) of Theorem \ref{BG} to the category of curved differential quadratic algebras which is ``one to one". This duality is summarized by the following theorem.

\begin{theorem}\label{NHKD}
The above correspondence defines an anti-isomorphism between the category of nonhomogeneous quadratic algebras satisfying Conditions $\mathrm{(i)}$ and $\mathrm{(ii)}$ of Theorem \ref{BG} and the category of curved differential quadratic algebras which induces an anti-isomorphism between the category of nonhomogeneous quadratic Koszul algebras and the category of curved differential quadratic Koszul algebras.
\end{theorem}

\noindent \underbar{Remark}. It is important to realize that an object $\fraca$ of the category of nonhomogeneous quadratic algebras as defined above is not only the filtered algebra $\fraca$ but it is $\fraca=A(E,P)$ defined as a quotient of a tensor algebra $T(E)$ as explained before. In particular the vector space $E$ is explicitely involved although its canonical image in the filtered algebra is a quotient which is generically not isomorphic to $E$; see for instance in Section \ref{aalg} the case where $\psi:E\otimes E\rightarrow E$ is a nonassociative product. The morphisms of this category have been defined accordingly. Of course if $\fraca$ has the PBW property $E$ is mapped isomorphically onto the corresponding subspace of the algebra $\fraca$. More generally it is worth noticing here that if the conditions (i) and (ii) of Theorem \ref{BG} are satisfied, then the canonical projection $p:T(E)\rightarrow \fraca$ is injective on $E$, i.e. it induces an isomorphism $E\simeq p(E)$, 
\cite{bra-gai:1996} (see \S 3.2 in \cite{bra-gai:1996}).

\section{Examples II}\label{ExII}

In this section we analyze some examples of nonhomogeneous quadratic Koszul algebras having examples of Section \ref{ExI} as homogeneous parts and we describe their dual (curved) differential quadratic algebras.
\\

\subsection*{ 1. Universal enveloping algebra $U(\fracg)$}~\\
Let $\fracg$ be a Lie algebra and let $U(\fracg)$ be its universal enveloping algebra. It is clear that $\fraca=U(\fracg)$ is a nonhomogeneous quadratic algebra with quadratic part $\cala$ which coincides with the symmetric algebra $S(\fracg)$, $\cala=S(\fracg)$. It follows that for the Koszul dual one has $\cala^!=\wedge \fracg^\ast$. We know (see in Section \ref{ExI}, Example 1) that $S(\fracg)$ is Koszul of global dimension $D=\dim(\fracg)$ and is in fact Koszul-Gorenstein with polynomial growth. Furthermore it is a classical result that $U(\fracg)$ has the PBW property (see also in Section \ref{Lie}). Thus $U(\fracg)$ is a Koszul nonhomogeneous quadratic algebra. In this case $\psi_0=0$, i.e. there are only linear nonquadratic terms, so $\wedge \fracg^\ast$ is a differential quadratic Koszul algebra (curvature zero).
\\

\subsection*{ 2. The case where $E$ is an associative algebra $A$} ~\\
Let $A$ be an associative algebra with product $(x,y)\mapsto xy$ and let $\fraca$ be the nonhomogeneous quadratic algebra defined by
\[
\fraca=T(A)/(\{x\otimes y-xy\vert x,y\in A\})
\]
as in Section \ref{aalg}. Then the quadratic part is given by
\[
\cala=T(A)/(A\otimes A)=A\oplus \mathbb K\bbbone
\]
i.e. $\cala_1=A$ as vector space and $\cala_n=0$ for $n\geq 2$. The quadratic algebra $\cala$ is Koszul but is of infinite global dimension (see in Section \ref{ExI}, Example 2). In view of Theorem \ref{ALG}, $\fraca$ has the PBW-property and therefore $\fraca$ is a Koszul nonhomogeneous quadratic algebra.\\ Since $\psi_0=0$, the Koszul dual $\cala^!$ of $\cala$ which is the tensor algebra over the dual vector space $A^\ast$ of $A$, i.e. $\cala^!=T(A^\ast)$, is a differential quadratic Koszul algebra, (see also Theorem \ref{ALG} b)).
\\

\subsection*{ 3. Calculus on the twisted $SU(2)$ group \cite{wor:1987b}} ~\\
Let $\fraca$ be the nonhomogeneous quadratic algebra generated by 3 elements $\nabla_0, \nabla_1, \nabla_2$ with relations

\begin{equation}
\left\{
\begin{array}{l}
\nu\nabla_2\nabla_0-\displaystyle{\frac{1}{\nu}}\nabla_0\nabla_2 = \nabla_1\\
\\
\nu^2\nabla_1\nabla_0-\displaystyle{\frac{1}{\nu^2}}\nabla_0\nabla_1 = (1+\nu^2)\nabla_0\\
\\
\nu^2\nabla_2\nabla_1-\displaystyle{\frac{1}{\nu^2}}\nabla_1\nabla_2 = (1+\nu^2)\nabla_2
\end{array}
\right.
\label{rWn}
\end{equation}
with $\nu\in \mathbb K\backslash \{0\}$. The quadratic part of $\fraca$ is the quadratic algebra $\cala$ of Example 3 in Section \ref{ExI} with relations given by (\ref{rWh}), (i.e. the above relations where the non quadratic terms of the right-hand sides have been set equal to 0). This quadratic algebra is Koszul of global dimension 3 and is Koszul-Gorenstein with polynomial growth. It is not hard to verify that $\fraca$ has the PBW property by using Theorem \ref{BG} and Proposition \ref{EPBW} with $\psi_0=0$ (since the nonquadratic terms are linear). It follows that the Koszul dual $\cala^!$ of $\cala$, which is as explained in Section \ref{ExI} the quadratic algebra generated by $\omega_0, \omega_1, \omega_2$ in degree 1 with relations (\ref{KdrWR}), is a differential quadratic Koszul algebra; its differential is the antiderivation $d$ given by
\begin{equation}
\left\{
\begin{array}{l}
d\omega_0=\nu^2(1+\nu^2)\omega_0\omega_1\\
d\omega_1=\nu \omega_0\omega_2\\
d\omega_2=\nu^2 (1+\nu^2)\omega_1\omega_2
\end{array}
\right.
\end{equation}
on the generators which implies $d^2=0$ on $\cala^!$. This is the left covariant differential calculus on the twisted $SU(2)$ group of \cite{wor:1987b}. More precisely it is a calculus on the twisted $SL(2)$ group and it becomes a calculus for the twisted $SU(2)$ group if one equips it with the involution 
\[
\nabla_1\mapsto \nabla^\ast_1=\nabla_1,\  \nabla_0\mapsto \nabla^\ast_0=-\nu^{-1}\nabla_2,\  \nabla_2\mapsto \nabla^\ast_2=-\nu\nabla_0
\]
and the corresponding one on the dual basis of the $\omega$'s.
\\

\subsection*{ 4. Clifford algebra and generalization}~\\
Given a symmetric bilinear form $g$ on $E$, the Clifford algebra of $g$
\[
\fraca=T(E)/(\{x\otimes y+y\otimes x-2g(x,y)\bbbone\vert x,y\in E\})
\]
gives another example of nonhomogeneous quadratic Koszul algebra. In this case the quadratic part $\cala$ of $\fraca$ is the exterior algebra
\[
\cala=\wedge E
\]
with Koszul dual $\cala^!=S(E^\ast)$ so it is of infinite global dimension. The algebra $\cala^!$ is a curved differential quadratic Koszul algebra with curvature $F=-2g\in S^2(E^\ast)$ but with trivial differential $d=0$.\\

More generally it is worth noticing that an antiderivation $d$ of $S(E^\ast)$ vanishes in degrees different from 1 since
\[
d(\alpha\beta)=d(\alpha)\beta-\alpha d(\beta)=-(d(\beta)\alpha-\beta d(\alpha))=-d(\beta\alpha)=-d(\alpha\beta)=0
\]
for any $\alpha, \beta\in E^\ast$. Furthermore
\[
d(\alpha)\beta=d(\beta)\alpha
\]
for any $\alpha, \beta\in E^\ast$ implies that there is a $\theta\in E^\ast$ such that one has
\[
d(\alpha) + \theta \alpha=0
\]
for any $\alpha\in E^\ast$ whenever $d$ is of degree 1. Using Theorem \ref{NHKD} it follows that the algebras of the form
\[
\fraca=T(E)/(\{ x\otimes y+ y\otimes x -\theta (x)y-\theta(y)x-2g(x,y)\bbbone\})
\]
with $\theta\in E^\ast$ and $g\in S^2(E^\ast)$ are all the nonhomogeneous quadratic Koszul algebras having quadratic part $\cala=\wedge E$.

\section{Lie prealgebras}\label{pseudoL}

Throughout the following, $\fraca=T(E)/(P)$ is a nonhomogeneous quadratic algebra with homogeneous part $\cala=T(E)/(R)$ and we assume that the conditions (i) and (ii) of Theorem \ref{BG} are satisfied with $\psi_1:R\rightarrow E$ and $\psi_0:R\rightarrow \mathbb K$ as before but now we take $\psi_0=0$. This latter assumption is equivalent to say that $\fraca$ is an {\sl augmented algebra}, its augmentation $\varepsilon:\fraca\rightarrow \mathbb K$ being induced by the projection of $T(E)$ onto the degree zero component ($\simeq \mathbb K$). In this case $F=0$ so $\cala^!$ is a graded differential algebra in view of (\ref{Curvd}) with differential $d$ induced by the transpose of $\psi_1$. An algebra $\fraca$ as above is referred to as a {\sl quadratic-linear algebra}, \cite{pol-pos:2005}. Notice that then $\fraca$ is an augmented filtered algebra and that the subspace of $F^1(\fraca)$ which is in the kernel of the augmentation, i.e. $p(E)=F^1(\fraca)\cap \ker(\varepsilon)$ is the image of $E$ under the canonical projection $p:T(E)\rightarrow \fraca$. In the case where $\fraca$ has the PBW property, $E$ can be identified with its image $p(E)$  and this is in fact already true here since we have assumed that the conditions (i) and (ii) of Theorem \ref{BG} are satisfied \cite{bra-gai:1996} (see in the remark at the end of Section \ref{dual}).  Theorem \ref{NHKD} has the following counterpart.

\begin{theorem}\label{QLKD}
The anti-isomorphism of Theorem \ref{NHKD} induces an anti-isomor\-phism between the category of quadratic-linear algebras and the category of differential quadratic algebras which restricts as an anti-isomorphism between the category of quadratic-linear Koszul algebras and the category of differential quadratic Koszul algebras.
\end{theorem}

\noindent \underbar{Example}. Let $\fracg$ be a Lie algebra then $\fraca=U(\fracg)$ is of the above type with $\cala=S(\fracg)$ so $\cala^!=\wedge \fracg^\ast$. 
As explained in Section \ref{Lie} this gives an example of the above general situation.
Furthermore the correspondence $\fracg\mapsto \wedge\fracg^\ast$ is an anti-isomorphism of categories between the category of Lie algebras and the category of free graded-commutative differential algebras generated in degree 1 (i.e. of the exterior algebras equipped with differentials).\\

This example shows that the class of augmented nonhomogeneous quadra\-tic algebras $\fraca$ considered here generalizes the class of the universal enveloping algebras of Lie algebra.  Notice that this class contains also the algebra $\fraca=T(A)/(\{x\otimes y-xy\vert x,y\in A\})$ associated to an associative algebra $A$ as explained in Section \ref{aalg}. This latter kind of quadratic-linear algebra is interesting but is clearly not a good analogue of universal enveloping algebras of Lie algebras. This is why we shall restrict attention to a smaller class of algebras where we impose the  Poincar\'e duality corresponding to the Gorenstein property and the assumption of Koszulity which is natural and always satisfied by the examples we have in mind. \\

Let $\fraca$ be a Koszul nonhomogeneous quadratic algebra which is quadratic-linear and such that its quadratic part $\cala$ is Koszul-Gorenstein (see in Section \ref{quad}) and let $E$ be the vector space of elements of degree 1 of $\cala$ identified with the corresponding generating subspace of $\fraca$, that is 
\[
E=F^1(\fraca)\cap \ker (\varepsilon)
\]
where $\varepsilon:\fraca\rightarrow \mathbb K$ is the augmentation of $\fraca$. Under these assumptions, $E$ will be refered to as a {\sl Lie prealgebra} and $\fraca$ will be refered to as {\sl the enveloping algebra of $E$}
. If, furthermore the quadratic part $\cala$ has polynomial growth (see remark b) in Section \ref{ExI}) we shall speak of a {\sl regular Lie prealgebra}.
Examples 1 and 3 of last section give examples of regular Lie prealgebras while Example 2 which is quadratic-linear and Koszul does not correspond to a Lie prealgebra.\\

We formalize the notion of a Lie prealgebra with the following definition.

\begin{definition}
A {\sl Lie prealgebra} is a triple $(E,R_E,\psi_E)$ where $E$ is a finite-dimensional vector space, $R_E$ is a vector subspace of $E\otimes E$ and where 
\[
\psi_E:R_E\rightarrow E
\]
 is a linear mapping such that the quadratic algebra
\[
\cala_E=T(E)/(R_E)
\]
is Koszul-Gorenstein and such that the quadratic-linear algebra
\[
\fraca_E=T(E)/(\{r+\psi_E(r)\vert r\in R_E\})
\]
has the PBW property (so is Koszul). When no confusion arises concerning $R_E$ and $\psi_E$ the Lie prealgebra will be simply denoted by $E$; the quadratic-linear algebra $\fraca_E$ is its {\sl enveloping algebra}. A {\sl morphism} $\alpha$ from $(E,R_E,\psi_E)$ to another Lie prealgebra $(E',R_{E'},\psi_{E'})$ is a linear mapping $\alpha:E\rightarrow E'$ such that $(\alpha\otimes \alpha) (R_E)\subset R_{E'}$ and $\alpha\circ \psi_E=\psi_E'\circ(\alpha\otimes \alpha)$.
\end{definition}

This defines the category of Lie prealgebras and $E\mapsto \fraca_E$ is a covariant functor from the category of Lie prealgebras to the category of quadratic-linear Koszul algebras.\\

We now wish to describe completely the duality of Theorem \ref{QLKD} restricted to Lie prealgebras.\\
Let us remind that a {\sl Frobenius algebra} is a finite-dimensional algebra $\cala$ such that one has the left $\cala$-module isomorphism $\cala\simeq {\cala^\ast}$ between $\cala$ and its dual vector space $\cala^\ast$ where the left $\cala$-module structure of $\cala^\ast$ is induced by the right $\cala$-module structure of $\cala$. This is equivalent to the existence of a nondegenerate bilinear form $(x,y)\mapsto B(x,y)$ on $\cala$ such that 
\[
B(xy,z)=B(x,yz)
\]
for $x,y,z\in \cala$. Then the linear mapping 
\[
x\mapsto B(\bullet, x)
\]
realizes the left $\cala$-module isomorphism $\cala\simeq \cala^\ast$. Concerning the graded connected case, one has the following result \cite{smi:1996}.

\begin{proposition}\label{GFA}
Let $\cala$ be a finite-dimensional graded connected algebra such that $\cala_D\not=0$ and $\cala_n=0$ for $n>D$. The following conditions $\mathrm{(i)}$ and $\mathrm{(ii)}$ are equivalent :\\
$\mathrm{(i)}$ $\cala$ is a Frobenius algebra,\\
$\mathrm{(ii)}$ $\dim(\cala_D)=1$ and $(x,y)\mapsto (xy)_D$ is nondegenerate,\\
where $(xy)_D$ is the component of $xy$ in $\cala_D$.\\
If $\cala$ is a graded connected Frobenius algebra, then there is an automorphism $\sigma\in \Aut(\cala)$ of graded algebra (i.e. of degree 0) such that
\[
xy=\sigma(y)x
\]
for any $x\in \cala_n$ and $y\in \cala_{D-n}$ with $D$ as above and $0\leq n\leq D$.
\end{proposition}

\noindent \underbar{Examples}. The prototype of a graded Frobenius algebra as above is the exterior algebra $\wedge E$ over a finite-dimensional vector space $E$ with $D=\dim(E)$.\\
More generally, the Koszul dual $\cala^!$ of the quadratic algebra $\cala$ of Example 4 in Section \ref{ExI}, that is the quadratic algebra generated by $D$ elements $\omega_\alpha$ ($\alpha\in \{1,\dots,D\}$) with relations given by (\ref{qwedgeD}) with the $q$'s satisfying (\ref{qrs}) is a graded Frobenius algebra as above. Indeed a basis of this algebra is given by the $\omega_{i_1}\dots \omega_{i_k}$ with $i_1<i_2<\dots<i_k$ and one verifies that Condition (ii) of Proposition \ref{GFA} is satisfied. In this case, the automorphism $\sigma$ of that proposition is given by
\begin{equation}
\sigma(\omega_\alpha)=(-1)^{D-1}\prod_{\beta\not=\alpha} q^{\alpha\beta} \omega_\alpha
\label{sigma}
\end{equation}
for $\alpha\in\{1,\dots,D\}$.\\
A particular case of the above algebra is the exterior algebra over a $D$-dimensional vector space $E$ which corresponds to $q^{\alpha\beta}=1$ for all $\alpha, \beta\in \{1,\dots,D\}$. Formula (\ref{sigma}) gives then
\[
\sigma(x)=(-1)^{(D-1)n}x
\]
for $x\in \wedge^nE$, $0\leq n\leq D$.\\
Another particular case of the above algebra with $D=3$ is the quadratic algebra $\cala^!$ with relations (\ref{KdrWR}) in (Example 3) Section \ref{ExI}. Formula (\ref{sigma}) reads
\[
\sigma(\omega_0)=\nu^6\omega_0,\ \ \sigma (\omega_1)=\omega_1,\ \  \sigma(\omega_2)=\nu^{-6}\omega_2
\]
for the automorphism $\sigma$.\\

In the following we shall consider differential quadratic (Koszul) algebras which are Frobenius and we shall call such an algebra a differential quadratic (Koszul) Frobenius algebra. One has the corresponding category. It turns out that the duality of Theorem \ref{QLKD} restricts as a duality between Lie prealgebras and differential quadratic Koszul Frobenius algebras.

\begin{theorem}\label{PLK}
Let $E$ be a Lie prealgebra then $\cala^!_E$ is a differential quadratic Koszul Frobenius algebra and this defines an anti-isomorphism between the category of Lie prealgebras and the category of differential quadratic Koszul Frobenius algebras.
\end{theorem}

This follows from the fact proved in \cite{smi:1996} that a quadratic Koszul algebra $\cala$ of finite global dimension is Koszul-Gorenstein if and only if its Koszul dual $\cala^!$ is Frobenius.

\section{Complexes and homologies}\label{coh}

In this section we define the generalization of the Chevalley-Eilenberg complexes for a Lie prealgebra and relate their (co)homologies to the corresponding Hochschild (co)homologies of its enveloping algebra generalizing thereby the relation between the Chevalley-Eilenberg (co)homologies of a Lie algebra $\fracg$ and the corresponding Hochschild (co)homologies of its universal enveloping algebra $U(\fracg)$. It turns out that several results of these developments are valid in the more general context of quadratic-linear algebras. Let us complete accordingly the definitions of last section.

\begin{definition}
A {\sl prealgebra} is a triple $(E,R_E,\psi_E)$ where $E$ is a finite dimensional vector space, $R_E$ is a vector subspace of $E\otimes E$ and where $\psi_E:R_E\rightarrow E$ is a linear mapping such that
\[
\fraca_E=T(E)/(\{r+\psi_E(r)\vert r\in R_E\})
\]
is a quadratic-linear algebra, that is such that $\psi_1=\psi_E$ and $\psi_0=0$ satisfy the conditions of Proposition \ref{EPBW}. When no confusion arises concerning $R_E$ and $\psi_E$ the prealgebra will be simply denoted by $E$; the quadratic linear algebra $\fraca_E$ is its {\sl enveloping algebra}. The morphisms of prealgebras are defined as in the case of Definition 1 for Lie prealgebras.
\end{definition}

One has therefore a category of prealgebras and the full subcategory of Lie prealgebras. Notice that the functor $E\mapsto \fraca_E$ is such that one has $E\subset \fraca_E$ canonically and that one can recover the prealgebra $(E, R_E,\psi_E)$ from the (augmented filtered) quadratic-linear algebra $\fraca_E$.\\

Given a prealgebra $(E,R_E,\psi_E)$ we let $\cala_E=A(E,R_E)$ be the quadratic part of its enveloping algebra $\fraca_E$. One can reformulate Theorem \ref{QLKD} in the following manner similar to Theorem  \ref{PLK}.

\begin{theorem} \label{PAK}
Let $E$ be a prealgebra then $\cala^!_E$ is a differential quadratic algebra and this defines an anti-isomorphism between the category of prealgebras and the category of differential quadratic algebras.
\end{theorem}

A {\sl representation} of the prealgebra $E=(E,R_E,\psi_E)$ also called a {\sl left $E$-module} is a vector space $V$ equipped with a linear mapping $\pi$ of $E$ into $\End(V)$ such that 
\begin{equation}
m(\pi\otimes \pi (r)) + \pi(\psi_E(r))=0
\label{reppa}
\end{equation}
for any $r\in R_E$ where $m:\End (V) \otimes \End(V)\rightarrow \End(V)$ denotes the product of endomorphisms. This is the same thing as a left $\fraca_E$-module since it extends as a representation of $T(E)$ which vanishes on $\{r+\psi_E(r)\vert r\in R_E\}$ in view of (\ref{reppa}).\\

One defines similarily a {\sl right representation} of $E$ or a {\sl right $E$-module} as corresponding to a right $\fraca_E$-module; for instance if $V$ is a left $E$-module, the dual vector space $V^\ast$ is canonically a right $E$-module.\\

Let $(\pi,V)$ be a representation of $E$ as above and let
\[
\delta^{(0)}_\pi:V\rightarrow V\otimes E^\ast
\]
be the linear mapping defined by
\[
\delta^{(0)}_\pi(v)(e)=\pi(e) v
\]
for $e\in E$ and $v\in V$. The linear mapping $\delta^{(0)}_\pi$ 
extends canonically as a right $\cala^!_E$-module endomorphism
\[
\delta_\pi:V \otimes \cala^!_E\rightarrow V\otimes \cala^!_E
\]
which is in fact an endomorphism of degree 1 of the free graded right $\cala^!_E$-module $V\otimes \cala^!_E$.\\
Let $d$ be the differential of $\cala^!_E$ and let us denote again by $d$ the linear mapping of degree 1 of $V\otimes \cala^!_E$ into itself defined by
\[
d(v\otimes w)=v\otimes d\omega
\]
for $v\in V$ and $\omega\in \cala^!_E$. One has the following result.

\begin{lemma}\label{DV}
The linear mapping of degree 1
\[
\delta_\pi+d:V \otimes \cala^!_E\rightarrow V \otimes \cala^!_E
\]
satisfies 
\[
(\delta_\pi+d)^2=0
\]
i.e. the graded vector space $V\otimes \cala^!_E$ endowed with $d_\pi=\delta_\pi+d$ is a cochain complex of vector spaces.
\end{lemma}

\noindent\underbar{Proof}. One has
\[
(\delta_\pi+d)v\otimes \omega = \delta_\pi(v)\omega+v\otimes d\omega
\]
for $v\in V$ and $\omega\in \cala^!_E$. This implies 
\[
\begin{array}{lll}
(\delta_\pi+d)^2 v\otimes \omega & = &
\delta^2_\pi(v)\omega+\delta_\pi(v)d\omega+d(\delta_\pi(v))\omega-\delta_\pi(v)d\omega+v\otimes d^2\omega\\
& = & (\delta^2_\pi(v)+d(\delta_\pi(v)))\omega.
\end{array}
\]
On the other hand equation (\ref{reppa}) is equivalent to
\begin{equation}
\delta^2_\pi(v)+d\delta_\pi(v)=0
\label{dreppa}
\end{equation}
for any $v\in V$.~$\square$\\

The cochain complexes $(V\otimes \cala^!_E,d_\pi)$ generalize the Chevalley-Eilenberg cochain complexes. Indeed if $E$ is a Lie algebra $\fracg$ then $\cala^!_E=\wedge \fracg^\ast$ and if $(\pi,V)$ is a representation of $\fracg$ then $(V \otimes \wedge \fracg^\ast,d_\pi)$ is a Chevalley-Eilenberg cochain complex and all these cochain complexes are of this form.\\

By duality one obtains a chain complex $(V^\ast\otimes \cala^{!\ast}_E, d^\ast_\pi)$ from the complex $(V\otimes \cala^{!}_E,d_\pi)$, in particular one has the chain complex $(\cala^{!\ast}_E,d^\ast)$ which is a differential coalgebra. More generally, given a right representation $(W,\pi)$ one has the chain complex $(W\otimes \cala^{!\ast}_E,d^\ast_\pi)$. These chain complexes generalize the Chevalley-Eilenberg chain complexes.\\

We now assume in the sequel of this section that $E=(E,R_E,\psi_E)$ is a Lie prealgebra, i.e. that the enveloping algebra $\fraca_E$ is Koszul and that its quadratic part $\cala_E$ is Koszul-Gorenstein. This is equivalent to say that $\cala^!_E$ is a differential quadratic Koszul Frobenius algebra in view of Theorem \ref{PLK}. In this case, one has
\begin{equation}
H^\bullet (V\otimes \cala^!_E)=\Ext^\bullet_{\fraca_E}(\mathbb K, V)
\label{HExt}
\end{equation}
for a representation $V$ of $E$ and
\begin{equation}
H_\bullet(W\otimes \cala^{!\ast}_E)=\text{Tor}_\bullet^{\fraca_E}(W,\mathbb K)
\label{HTor}
\end{equation}
for a right representation of $E$ \cite{pol-pos:2005}, where $H^\bullet (V\otimes \cala^!_E)$ is the cohomology of $(V\otimes \cala^!_E,d_\pi)$ and $H_\bullet(W\otimes \cala^{!\ast}_E,d^\ast_\pi)$ is the homology of $(W\otimes \cala^{!\ast}_E,d^\ast_\pi)$. This generalizes what happens in the case where $E$ is a Lie-algebra $\fracg$, \cite{car-eil:1973}, \cite{lod:1992}, \cite{wei:1994}.\\

Recalling that for an algebra $\fraca$ and a $(\fraca,\fraca)$-bimodule $\calm$,  that is a $\fraca^e$-module $\calm$ (with $\fraca^e=\fraca\otimes \fraca^{opp}$), one has
\[
H^\bullet (\fraca,\calm)=\Ext^\bullet_{\fraca^e}(\fraca,\calm)
\]
for the Hochschild cohomology $H^\bullet(\fraca,\calm)$ of $\fraca$ and 
\[
H_\bullet (\fraca,\calm)=\text{Tor}_\bullet^{\fraca^e}(\calm,\fraca)
\]
for the Hochschild homology $H_\bullet(\fraca,\calm)$ of $\fraca$ with value in $\calm$, and remembering that $\mathbb K$ can be considered as a left or as a right $\fraca_E$-module by using the augmentation $\varepsilon:\fraca_E\rightarrow \mathbb K$, one deduces from (\ref{HExt}) and (\ref{HTor}) relations between the above (co)homologies (of $E$) and the Hochschild (co)homologies of the enveloping algebra $\fraca_E$ of $E$ which generalizes the relations between the Chevalley-Eilenberg (co)homologies of a Lie algebra $\fracg$ and the Hochschild (co)homologies of its universal enveloping algebra $U(\fracg)$.

\section{Further prospect (tentative conclusion)}

It is natural to expect that the differential calculi on quantum groups fall in the framework described above namely the Koszul duality of Positselski starting from Lie prealgebras.\\
In fact, it is usually not difficult to show that the quantum tangent space $E$ (see e.g. in \cite{kli-sch:1997}, Chapter 14 for the definition) is a prealgebra in the sense of Definition 2 and that the corresponding differential quadratic algebra $\cala^!_E$ is a Frobenius algebra. Indeed in most cases $\cala^!_E$ is explicitely given as a differential quadratic algebra which is by inspection a Frobenius algebra. However it is generally more involved to prove the Koszulity of the enveloping algebra $\fraca_E$ that is here the Koszulity of its quadratic part $\cala_E$.\\
In other words, it is generally easy to show that the quantum tangent spaces are {\sl weak Lie prealgebras} in the sense of the following definition.

\begin{definition}
A weak Lie prealgebra is a prealgebra $E$ such that the differential quadratic algebra $\cala^!_E$ is a Frobenius algebra.
\end{definition}

Given such a weak Lie prealgebra $E$, it is a Lie prealgebra if and only if $\cala_E$ is Koszul, or equivalently if and only if $\cala^!_E$ is Koszul.\\

For instance we have shown that relations (\ref{rWn}) define a Lie prealgebra; it turns out that this is the quantum tangent space corresponding to the left covariant differential calculus on the twisted $SL(2)$ group of \cite{wor:1987b}. On the other hand one shows easily that the quantum tangent spaces corresponding to the 4D$_\pm$ bicovariant differential calculi on the same twisted (quantum) $SL(2)$ group are weak Lie prealgebras. We conjecture that they are in fact Lie prealgebras but to show this one has to prove the Koszulity of the corresponding algebras.\\

It is worth noticing that the Koszul property is a very nice and desirable property which here is furthermore natural in view of the duality described in Theorem \ref{PLK} as well as in view of the (co)homological results described in last section. In contrast, Definition 3 is somehow artificial and corresponds to intermediate properties. Thus an important task remains to show that the usual quantum tangent spaces corresponding to the differential calculi on quantum groups are Lie prealgebras.\\

Finally we mention that we are aware that there are many very interesting related papers on the so-called quantum Lie algebras. However, to the best of our knowledge the description given in the present paper of the basic structures which are relevant is new.

\bibliographystyle{amsalpha}

\begin{thebibliography}{A}

\bibitem{art-sch:1987}
M.~Artin and W.F. Schelter.
\newblock Graded algebras of global dimension 3.
\newblock {\em Adv. Math.}, 66:171--216, 1987.

\bibitem{art-tat-vdb:1991}
M.~Artin, J.~Tate, and M.~Van~den Bergh.
\newblock Modules over regular algebras of dimension 3.
\newblock {\em Invent. Math.}, 106:335--388, 1991.


\bibitem{ber:2005}
R.~Berger.
\newblock Dimension de {H}ochschild des alg{\`e}bres gradu{\'e}es.
\newblock {\em C.R. Acad.Sci. Paris, Ser. I}, 341:597--600, 2005.

\bibitem{bra-gai:1996}
A.~Braverman and D.~Gaitsgory.
\newblock Poincar{\'e}-{B}irkhoff-{W}itt theorem for quadratic algebras of
  {K}oszul type.
\newblock {\em J. Algebra}, 181:315--328, 1996.

\bibitem{car:1958}
H.~Cartan.
\newblock Homologie et cohomologie d'une alg{\`e}bre gradu{\'e}e.
\newblock {\em S{\'e}minaire {H}enri {C}artan}, 11(2):1--20, 1958.

\bibitem{car-eil:1973}
H.~Cartan and S.~Eilenberg.
\newblock {\em Homological algebra}.
\newblock Princeton University Press, 1973.

\bibitem{mdv:2001}
M.~Dubois-Violette.
\newblock Lectures on graded differential algebras and noncommutative geometry.
\newblock In Y.~Maeda et. al. editors, {\em Noncommutative Differential
  Geometry and Its Applications to Physics}, pages 245--306. Shonan, Japan,
  1999, Kluwer Academic Publishers, 2001.

\bibitem{flo:2006}
G.~Fl{\o}ystad.
\newblock Koszul duality and equivalences of categories.
\newblock {\em Trans. Amer. Math. Soc.}, 358:2373--2398, 2006.

\bibitem{gur:1990}
D.I. Gurevich.
\newblock Algebraic aspects of the quantum {Y}ang-{B}axter equation.
\newblock {\em Algebra i Analiz (Transl. in Leningrad Math. J. 2 (1991)
  801-828)}, 2:119--148, 1990.

\bibitem{kli-sch:1997}
A.~Klimyk and K.~Schm{\"u}dgen.
\newblock {\em Quantum groups and their representations}.
\newblock Springer, 1997.

\bibitem{lod:1992}
J.L. Loday.
\newblock {\em Cyclic homology}.
\newblock Springer Verlag, New York, 1992.

\bibitem{man:1987}
Yu.~I. Manin.
\newblock Some remarks on {K}oszul algebras and quantum groups.
\newblock {\em Ann. Inst. Fourier, Grenoble}, 37:191--205, 1987.

\bibitem{man:1988}
Yu.~I. Manin.
\newblock {\em Quantum groups and non-commutative geometry}.
\newblock CRM Universit{\'e} de Montr{\'e}al, 1988.

\bibitem{pol-pos:2005}
A.~Polishchuk and L.~Positselski.
\newblock {\em Quadratic algebras}, volume~37 of {\em University {L}ecture
  {S}eries}.
\newblock Amer. {M}ath. {S}oc., Providence, RI., 2005.

\bibitem{pos:1993}
L.~Positselski.
\newblock Nonhomogeneous quadratic duality and curvature.
\newblock {\em Func. Anal. Appl.}, 27:197--204, 1993.

\bibitem{pot:2006}
A.~Pottier.
\newblock Stabilit\'e de la propri\'et\'e de {K}oszul pour les alg\`ebres
  homog\`enes vis-\`a-vis du produit semi-crois\'e.
\newblock {\em C.R. Acad. Sci. Paris, S{\'e}rie {I}}, 343:161--164, 2006.

\bibitem{smi:1996}
S.P. Smith.
\newblock Some finite dimensional algebras related to elliptic curves.
\newblock {\em CMS Conf. Proc. Proc.}, 19:315--348, 1996.

\bibitem{wam:1993}
M.~Wambst.
\newblock Complexes de {Koszul} quantiques.
\newblock {\em Ann. Inst. Fourier, Grenoble}, 43:1083--1156, 1993.

\bibitem{wei:1994}
C.A. Weibel.
\newblock {\em An introduction to homological algebra}.
\newblock Cambridge University Press, 1994.

\bibitem{wor:1987b}
S.L. Woronowicz.
\newblock Twisted {SU}(2) group. an example of noncommutative differential
  calculus.
\newblock {\em Publ. RIMS, Kyoto Univ.}, 23:117--181, 1987.

\bibitem{wor:1989}
S.L. Woronowicz.
\newblock Differential calculus on compact matrix pseudogroups (quantum
  groups).
\newblock {\em Commun. Math. Phys.}, 122:129--170, 1989.


\end{thebibliography}

\end{document}